%% file: main.tex
\documentclass[10pt, journal]{IEEEtran}
\IEEEoverridecommandlockouts

\IEEEoverridecommandlockouts

\usepackage{amsmath,amssymb,amsfonts,amsthm}
\usepackage{algorithm,algorithmicx}
\usepackage[noend]{algpseudocode} 
\usepackage{graphicx}
\usepackage{textcomp}
\usepackage{epstopdf}
\usepackage{xcolor}
\usepackage{cite}
\usepackage{latexsym,todonotes}
\usepackage{graphicx}
\usepackage{epstopdf}
\usepackage{float}
\usepackage{subfigure}
\usepackage{multirow}
\usepackage{txfonts}
\usepackage{mathptmx}
\usepackage{mathrsfs}
\usepackage{exscale}%
\usepackage{relsize}%

\usepackage{tikz,pgfplots}
\usepackage{flowchart}
\usepackage{tabularx}%
\usepackage{textcomp}
\usepackage{booktabs}%
\usepackage[switch]{lineno}%
\usepackage{bbding}%
\usepackage[justification=centering]{caption}
\usepackage[bookmarks=false]{hyperref}
\def\BibTeX{{\rm B\kern-.05em{\sc i\kern-.025em b}\kern-.08em
    T\kern-.1667em\lower.7ex\hbox{E}\kern-.125emX}}

\begin{document}

\title{A Ring Topology-based Communication-Efficient Scheme for D2D Wireless Federated Learning}

\author{\IEEEauthorblockN{ Zimu Xu, Wei Tian, Yingxin Liu, Wanjun Ning, and Jingjin Wu, \emph{Member, IEEE}}

\thanks{The authors are with the Department of Statistics and Data Science, BNU-HKBU United International College, Zhuhai, Guangdong, 519087, P. R. China. J. Wu is also with the Guangdong Provincial Key Laboratory of Interdisciplinary Research and Application for Data
Science, Guangdong, 519087, P. R. China. (E-mail: \href{mailto:r130202603@mail.uic.edu.cn}{r130202603@mail.uic.edu.cn}; \href{mailto:s230202702@mail.uic.edu.cn}{s230202702@mail.uic.edu.cn}; \href{mailto:p930026088@mail.uic.edu.cn}{p930026088@mail.uic.edu.cn}; \href{mailto:janinening@ieee.org}{janinening@ieee.org};  \href{mailto:jj.wu@ieee.org}{jj.wu@ieee.org}). Corresponding author: J. Wu. 

This work is supported by Zhuhai Basic and Applied Basic Research Foundation Grant ZH22017003200018PWC, the Guangdong Provincial Key Laboratory of 
Interdisciplinary Research and Application for Data Science, 
Project code 2022B1212010006, and Guangdong Higher Education Upgrading
Plan (2021-2025) UIC R0400001-22.
}}

\maketitle

\input{paper}

\bibliographystyle{IEEEtran}
\bibliography{IEEEexample}
\vspace{12pt}

\end{document}

%% file: paper.tex
\begin{abstract}
Federated learning (FL) is an emerging technique aiming at improving communication efficiency in distributed networks, where many clients often request to transmit their calculated parameters to an FL server simultaneously. 
However, in wireless networks, the above mechanism may lead to prolonged transmission time due to unreliable wireless transmission and limited bandwidth. This paper proposes a communication scheme to minimize the uplink transmission time for FL in wireless networks. The proposed approach consists of two major elements, namely a modified Ring All-reduce (MRAR) architecture that integrates D2D wireless communications to facilitate the communication process in FL, and a modified Ant Colony Optimization algorithm to identify the optimal composition of the MRAR architecture. Numerical results show that our proposed approach is robust and can significantly reduce the  transmission time compared to the conventional star topology. Notably, the reduction in uplink transmission time compared to baseline policies can be substantial in scenarios applicable to large-scale FL, where client devices are densely distributed.   
\end{abstract}


\section{Introduction}

Thanks to the rapid growth in the number and popularity of user devices (UDs) capable of collecting and processing data, artificial intelligence (AI) techniques are more commonly applied in many aspects of daily life. Federated learning (FL)~\cite{li2020federated} is a class of distributed machine learning techniques that keeps all actual training data on local UDs while only exchanging specific model parameters with each other. By decoupling the training process to local UDs, FL enables collaborative learning of AI models without storing data in the server. As a result, the privacy of UDs is better preserved than traditional centralized mechanisms. This advantage in turns motivates more users to participate in FL and leads to a more accurate global model.

A base station (BS) can function as an FL server in wireless networks, while the UDs are considered clients. One obstacle to the wider application of FL in wireless networks is that it requires multiple rounds of model parameters exchange between the FL server and clients until the FL process converges to an accurate global model~\cite{li2020federated}. 
In every training round, the BS updates the global model by aggregating the local parameters sent by UDs generated from their respective local models. Then, the BS broadcasts the updated global model to the UDs~\cite{8917724}. Most UDs participating in the training round will transmit within a short timeframe for synchronization purposes. This leads to a substantial strain on the uplink, especially in urban scenarios where a single BS serves a high volumes of UDs for the FL process. As a result, the uplink transmission latency would become unacceptable, and the overall FL performance would be negatively impacted. 


There are two straightforward directions to tackle the above-mentioned problem in existing studies. One direction, referred as gradient compression, aims to reduce the data size of transmitted parameters to improve communication efficiency~\cite{9177084}
Another direction is to reduce the frequency of communication by restricting the number of training rounds such that fewer communications between the BS and UDs would be required~\cite{liu2020client}. However, both approaches would sacrifice the accuracy of the final model to some extent. 

It has been investigated in more recent studies that D2D communications can be an attractive alternative to improve the communication efficiency for FL in wireless networks, as it can effectively avoid the communication bottleneck problem when compared to the traditional server-client (star) topology. Most approaches aiming at maximizing the efficiency of wireless D2D communications involve optimal resource allocation or power control~\cite{7051236}. One study that focused on the training algorithm is~\cite{xing2021federated}, which proposed to implement the decentralized stochastic gradient descent (DSGD) in an FL D2D wireless network. 
However, DSGD generally requires more communication rounds to achieve the same level of training accuracy as in conventional algorithms such as FedSGD, FedAvg and FedProx, due to its distributed implementation.

We take a different perspective in this paper, by proposing a topology-based approach that incorporates a modified version of recently proposed Ring All-reduce (RAR) architecture~\cite{9796785} into FL in orthogonal frequency division multiple access (OFDMA) wireless networks. 
Unlike the approach in~\cite{xing2021federated}, our approach does not alter FL algorithms and thus does not incur extra communication rounds to converge. Instead, we focus on the parameter aggregation process in each training round to minimize the uplink transmission time, by identifying the most efficient ring topology for aggregating the parameters. 


Another benefit of our  modified RAR (MRAR)-based approach is enhanced resilience against link failures. In case of link failures, only a small amount of additional data is required for re-transmission to maintain the aggregation accuracy under the MRAR implementation. This preserves the advantage in communication efficiency of MRAR over the star topology.

Our main contributions are summarized as follows,
\begin{itemize}
    \item We incorporate the MRAR architecture to address a potential bottleneck problem caused by simultaneous transmissions that may hinder the large-scale application of FL in OFDMA wireless networks. Note that we are not proposing a new FL algorithm, but a robust and efficient communication scheme to aggregate parameters under conventional FL algorithms. Therefore, our proposed approach can potentially improve the overall efficiency and performance of various existing FL algorithms.
    \item We analytically prove that, the expected uplink transmission time by our proposed MRAR approach increases sublinearly with the density of UDs. This is a significant improvement over the conventional star topology with a linear increasing rate.
    \item We propose an ACO-based algorithm to overcome the computational complexity in identifying the optimal composition of the MRAR, which in turns gives the shortest transmission time amongst all possible ring topologies.
    \item We demonstrate by simulation results that our proposed method significantly outperforms  existing communication schemes in terms of transmission efficiency and robustness. Particularly, in scenarios where the UDs are densely distributed, the reduction in transmission time can substantial. Also, within the MRAR framework, the ACO-based algorithm achieves further reduction in transmission time over the greedy algorithm. 
\end{itemize}





\section{Preliminaries and System Model}

\subsection{Parameter aggregation in federated learning}
We consider that a set of $\mathcal{V} = \{1,…,K\}$ of $K$ UDs are participating in an FL training round. In mainstream FL implementations such as FedSGD, FedAvg, FedProx and federated learning with quantization (FedPAQ)~\cite{reisizadeh2020fedpaq}, UDs do not share their own data, but only exchange and aggregate trained parameters with each other. At the end of each training round, UDs take the aggregated parameters from the global model as the initial parameters in the next round. 

The following equation represents a common form of final aggregation and update of global parameters,
\begin{small}
\begin{equation}
     w_{\text{global}} = \sum_{k \in \mathcal{V}} \frac{|D_k|}{|D|} w_k, \label{global}
\end{equation}
\end{small}\ignorespacesafterend 
where $|D_k|$ is the size of locally trained data at UD $k$, $|D|$ is the total size of trained data from all UDs, $w_{\text{global}}$ is the aggregated parameters or soft labels, and $w_k$ is locally trained parameters or soft labels at UD $k$.

We consider the uplink of an OFDMA wireless network. When FL is implemented in such networks with a traditional star topology, 
UDs synchronously send their locally trained parameters to the BS at the end of each training round. The BS then updates the global model as in~\eqref{global} before broadcasting the updated model to UDs. The UDs initialize their local model with the global model at the beginning of the next round.  


\begin{figure*}[htbp]
\centering
     \includegraphics[scale=0.54]{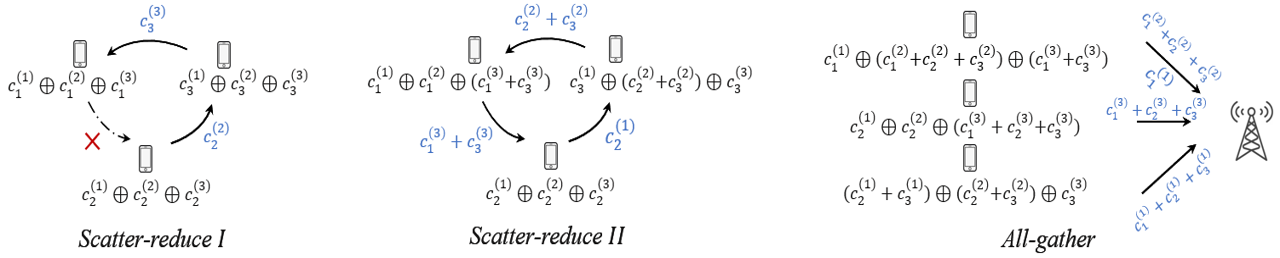}
     \vspace{0.1cm}
    \caption{\textit{An illustration of MRAR implementation with 3 UDs with a single-link failure}}
    \vspace{0.2cm}
    \label{interrupt}
\end{figure*}

\subsection{The Modified Ring All-reduce scheme}
We now introduce the MRAR scheme and show how it implements the aggregation task in FL in a different way. In MRAR, 
UDs form a logical ring where each UD has two neighbors. A UD will only send data to the ``next" neighbor in the clockwise direction and receive data from the ``last" neighbor in the counter-clockwise direction.


Two key steps to implement the MRAR are \textit{Scatter-reduce} and \textit{All-gather}. Assume that $K$ UDs participate in the current training round, each UD will equally divide its weighted local model parameter (i.e. $w_k$ in \eqref{global}) into $K$ chunks, that is, 
\begin{small}
\begin{equation}
   \frac{|D_k|}{|D|}w_k = \left(c_k^{(1)}\oplus c_k^{(2)}  \oplus ... \oplus c_k^{(K)} \right), 
\end{equation}
\end{small}\ignorespacesafterend
where $\oplus$ is the concatenation operation and $c_k^{(j)} = \frac{|D_k|}{|D|}w_k^{(j)}$ represents the $j$-th chunk for UD $k$.

If we label an arbitrary UD as UD $1$, and number the other UDs in ascending order according to the transmission direction, then the two steps in the $n$-th round can be described as follows,
\begin{itemize}

\item  \textit{Scatter-reduce}: UD $k$ sends the $\left[(k-n+1)\%K\right]$-th chunk of its accumulated weighted local model to its next neighbor, and receives its last neighbor's $\left[(k-n)\%K\right]$-th chunk. Here, $\%$ represents the modular operator. UD $k$ then accumulates the received chunk to obtain
\begin{small}
    \begin{equation}
        \underbrace{\sum_{j=k-n}^k c_{j\%K}^{((k-n)\%K)}}_{\text{To be sent in the next round}} = \underbrace{\sum_{j=k-n}^{k-1} c_{j\% K}^{((k-n)\%K)}}_{\text{The received chunk}} + c_k^{((k-n)\%K)}
    \end{equation}
\end{small}


\item \textit{All-gather}: UD $k$ sends the $\left[(k-1)\%K\right]$-th chunk to the BS, which then aggregates the global model by splicing all data chunks together. That is,
\begin{small}
\begin{equation}
    w_{\text{global}} = (\sum_{k=1}^K c_k^{(1)} \oplus \sum_{k=1}^K c_k^{(2)} \oplus ... \oplus \sum_{k=1}^K c_k^{(K)}).   
\end{equation}
\end{small}

If the link of UD $k$ to UD $k+1$ is down at the $n$th round, UD $k$ also sends the $\left[(k-n)\%K\right]$-th chunk to the BS, for recovering lost information due to previous link failure(s).
\end{itemize}

The entire implementation consists of $N-1$ steps of scatter-reduce followed by one all-gather. By the end of the process, the BS broadcast the global model to UDs for subsequent training rounds. Fig~\ref{interrupt} illustrates an example of the MRAR implementation with 3 participating UDs with a single link failure. In particular, Fig.~\ref{interrupt} also demonstrates that, at most one additional chunk  is required to be transmitted for a single link failure to obtain the same aggregation result at the BS. In this case, the additional chunk is $c_1^{(1)}$. 




\subsection{Communication model}
For uplink OFDMA systems, it is reasonable to consider the Signal-to-Noise Ratio (SNR) to  measure the quality of transmissions. Consider a transmission from UD $i$ to UD $j$ (we consider the BS as a special UD with index $0$ hereafter), the uplink SNR is $\displaystyle \text{SNR}_{i,j}=(d_{ij}^{-\alpha}p_i)/N_0,$ where $\alpha$ is the path loss exponent, $d_{i,j}$ is the distance between transmitter $i$ and receiver $j$, $p_i$ is the transmission power of UD $i$, and $N_0$ represents the noise power. Then, the data rate $R_{i, j}$ of transmitter $i$ to receiver $j$ can be represent by
$\displaystyle 
R_{i, j}=B_i\log(1+\text{SNR}_{i, j}),$ where $B_i$ is the allocated bandwidth for UD $i$. 

In terms of OFDMA, if we consider a worst-case scenario where all UDs request to transmit simultaneously, the allocated bandwidth should be constrained by 
$\sum_{i \in \mathcal{V}} B_i \leqslant B $, where $B$ is the total bandwidth allocated for FL.  

\subsection{Analysis of transmission time}

Recall the procedure of the traditional federated learning scheme, the aggregation (i.e., \eqref{global}) can only be performed by the BS after all UDs have finished transmissions. The uplink transmission time of each round is then 
\begin{small}
\begin{equation} 
T_{\text{star}} = \frac{M}{\min_{i \in \mathcal{V}}{R_{i, 0}}},
\end{equation}
\end{small}\ignorespacesafterend 
where $M$ is the model size.

In MRAR scheme, each UD will transmit $N-1$ chunks in scatter-reduce and $1$ chunk in all-gather. Since each chunk contains $1/N$ of model parameters, the transmission time of each round is
\begin{small}

\begin{equation} 
    T_{\text{MRAR}} = \underbrace{\frac{(N-1)M}{N\min_{i \in \mathcal{V}}{R_{i, r(i)}}}}_{T_{\text{SR}}}  + \underbrace{\max_{i \in \mathcal{V}} \frac{(\mathcal{I}_i+1)M}{NR_{i,0}}}_{T_{\text{AG}}},
\end{equation}\label{TMRAR}
\end{small}\ignorespacesafterend 
where $r(i)$ represents the index of the next neighbour of UD $i$ in the ring. $\mathcal{I}_i$ is the number of additional trunks transmitted to guarantee the aggregation accuracy, which is equal to the number of link failures between UD $i$ and $r(i)$ that occurred in scatter-reduce in the same training round. The two terms in~\eqref{TMRAR}, $T_{\text{SR}}$ and $T_{\text{AG}}$, refer to the transmission time in the scatter-reduce and all-gather steps, respectively.



Next, we show that the expected transmission time under the traditional star topology is approximately linearly correlated with the density of UDs participating in FL in the area, given that the distribution of UDs follows a spatial Poisson Point Process (PPP). In FL, the transmission time of each round is  determined by the slowest transmission. Therefore, we may optimally allocate the bandwidth to balance the transmission rate of all connections and thus 
minimize the transmission time. 
 If the UDs are distributed in a circle with radius $R$, with the BS at the centre, and the minimum distance between UD and BS is set to $1$, 
the expectation of transmission time is
\begin{small}
\begin{equation} 
    \begin{aligned}
       \mathbb{E}\left[T_{\text{star}}^*\right] &= \frac{ M}{B}\mathbb{E}\left[\sum_{k}\frac{1}{\log(1 + \text{SNR}_{k, 0})} \right] \\
        &= \frac{M}{B} \sum_{n = 1}^{\infty}  \frac{n[\lambda \pi R^2]^n e^{-\lambda \pi R^2}}{n!} \mathbb{E}\left[ \frac{1}{\log(1 + \text{SNR}_{k, 0})}\right] \\
        &= \frac{\lambda \pi R M}{B} \frac{2}{R^2} \int_{x=0}^{R} \frac{x dx}{\log (1+\frac{p x^{-\alpha}}{N_0})},
    \end{aligned} \label{ET}
\end{equation} 
\end{small}\ignorespacesafterend 
where $\lambda$ is the intensity of PPP, $p$ is the transmission power of each UD.
It can be inferred from~\eqref{ET} that the lower bound of $\mathbb{E}\left[T^*_{\text{star}} \right]$ also increases linearly with $\lambda$. 

We now show that the expected transmission time in MRAR increases sublinearly with $\lambda$. First, the transmission time in scatter-reduce through the optimal ring is bounded above by that through a greedy ring, which is constructed by the following procedures:
\begin{enumerate}
    \item Add an arbitrarily selected UD to an empty chain.
    \item Among UDs not yet in the chain, select the one with the shortest distance (and thus the highest SNR) with the last added UD. Add the selected UD to the chain. Repeat until all UDs are added.
    \item Connect the last UD to the first UD to form a ring. 
    
\end{enumerate} 

We derive the expected transmission time for the greedy ring, which is the upper bound for scatter-reduce as, 
\begin{small}
\begin{equation} 
    \begin{aligned}
        \mathbb{E}\left[T_{\text{SR}}^*\right] 
        & < \frac{M}{B}\mathbb{E}\left[\sum_k\frac{1}{\log(1+\text{SNR}_{k,r(k)})} + \frac{1}{\log(1+\frac{p(2R)^{-\alpha}}{N_0})}\right] \\ 
        &= \frac{M}{B}\underbrace{\mathbb{E}\left[\sum_k\frac{1}{\log(1+\text{SNR}_{k,r(k)})} \right]}_{C}  + \underbrace{\frac{M}{B\log(1+\frac{p(2R)^{-\alpha}}{N_0})}}_{b},
    \end{aligned} \label{UB}
\end{equation}
\end{small}\ignorespacesafterend 
where $\log(1+p(2R)^{-\alpha}/N_0)^{-1}$ represents the maximum transmission time of the link formed by step 3) in the greedy ring as $2R$ is the maximum distance between two arbitrary UDs. For \eqref{UB}, we further analyze the upper bound of $C$ as $b$ is a constant. Since the nodes generated by PPP are independent, $C$ can be transformed to
\begin{small}
\begin{equation}
    \begin{aligned}
        C = \mathbb{E} \left[\sum_k\int_{0}^{R}P'_k(x)\frac{Mdx}{\log(1 + \frac{pd_k^{-\alpha}}{N_0})}\right],
    \end{aligned} \label{chain}
\end{equation}
\end{small}\ignorespacesafterend 
where $P_k(x)$ is the probability that both UD $k$ and $r(k)$ are within a circle of radius $R$, and the distance between them is less than $x$. 
By the nature of PPP, we obtain
\begin{small}
\begin{equation}
    \begin{aligned}
          C &< \sum_{n=1}^{\infty} \frac{[\lambda \pi R^2]^n e^{-\lambda \pi R^2}}{n!} \mathbb{E} \left[\sum_{k=1}^{n}\frac{e^{-\lambda \pi x^2}Mdx}{\log(1+\frac{px^{-\alpha}}{N_0})}\right]\\ 
          &= \sum_{n=2}^{\infty} \frac{n[\lambda \pi R^2]^n e^{-\lambda \pi R^2}}{n!}   \int_{0}^{R} \frac{e^{-\lambda \pi x^2}dx}{\log(1 + \frac{px^{-\alpha}}{N_0})}
       \\&= \lambda \pi R^2  \int_{0}^{R} \frac{e^{-\lambda \pi x^2}dx}{\log(1 + \frac{px^{-\alpha}}{N_0})}.
    \end{aligned} \label{EC}
\end{equation} 
\end{small}

Therefore, the expected transmission time for scatter-reduce sublinearly increases with $\lambda$. For all-reduce, the analysis follows similarly from that for the star topology, that is
\begin{small}
\begin{equation}
    \mathbb{E}\left[T_{\text{AG}}^*\right] = \frac{ (\mathbb{E}\left[\sum_k\mathcal{I}_k\right]+1)M}{B} \int_{x=0}^{R} \frac{x dx}{\log (1+\frac{p x^{-\alpha}}{N_0})}. \label{EAG}
\end{equation}
\end{small}

In D2D-enabled wireless FL, link failures are rather infrequent as the data size per transmission is relatively small (only parameters), and the SNR is generally high due to shorter transmission distances. Therefore,  $\mathbb{E}\left[\sum_k\mathcal{I}_k\right]$ is close to $0$, and  $\mathbb{E}\left[T_{\text{AG}}^*\right]$ will remain approximately constant as $\lambda$ changes. Combining~\eqref{EC} and~\eqref{EAG}, we conclude that the expected total transmission time of MRAR increases sublinearly with $\lambda$.

\section{Optimal Formation of the MRAR Architecture}
In this section, we demonstrate the algorithm to identify the best ring to minimize the total transmission time, i.e. to find the next neighbor $r(k)$ for each UD $k$ such that 
\begin{small}
\begin{equation} 
    \begin{aligned}
        \min &\quad T_{\text{SR}}^* \\
        \text{s.t.} \quad\quad & r(i) \neq r(j), \quad  i, j\in \mathcal{V}, i \neq j \\ & r^{(n)}(k) = k, 
        \quad k \in \mathcal{V}. 
    \end{aligned} \label{Problem}
\end{equation}
\end{small}

Since \eqref{Problem} has the same constraints as the TSP and is also NP-hard, we apply the ACO, which is an efficient and the most commonly used method for solving the TSP~\cite{dorigo2006ant}, to solve this problem. The completed ACO algorithm is presented as Algorithm~\ref{algorithm1}, while the key steps of the algorithm are summarized as follows,

\begin{enumerate}
    \item Initialize $a$ ants on each UD.
    \item Each ant starts travelling from the current UD and constructs a ring by repeatedly applying the state transition rule,
\begin{small}
\begin{equation} 
    \begin{aligned}
        P_{i, j} = \frac{h_{i,j}^\beta R_{i,j}^\gamma}{\sum_{k \in \mathcal{V}} h_{i,k}^\beta R_{i,k}^\gamma},
    \end{aligned} \label{prob}
\end{equation}  
\end{small}\ignorespacesafterend   
where $\mathcal{V}$ contains UDs that have not been visited.

    \item Each ant calculates the total transmission time corresponding to the ring, and updates the pheromone by 
    \begin{small}
    \begin{equation} 
        \begin{aligned}
            h_{i, j} = \rho(h_{i, j} + \triangle h_{i, j}) + (1-\rho) \frac{1}{T_{\text{SR}^*}^*},
        \end{aligned} \label{update}
    \end{equation}
    \end{small}\ignorespacesafterend 
    where $\text{ring}^*$ is the so-far best ring that achieves minimum transmission time, and $\triangle h_{i, j}$ is the sum of $\frac{1}{T_{\text{SR}}^*}$ of all ants. 

    \item Stop the algorithm if it reaches the maximum number of iterations $t$, else go to step 1).
    
\end{enumerate}

\renewcommand{\algorithmicrequire}{\textbf{Input:}}  
\renewcommand{\algorithmicensure}{\textbf{Output:}} 
\vspace{-0.05cm}
\begin{algorithm}[h]  
  \caption{Ant colony optimization}  
  \begin{algorithmic}[1]  
    \Require
        $K$: the number of UDs;
    \Statex \quad \ \
        $d$: distance matrix;
    \Ensure
        $r^*:\mathcal{V} \rightarrow \mathcal{V}$,
    \Statex \quad \quad \
        i.e. the optimal connection scheme of the ring; 
    \State initialize $\textbf{h} = \textbf{1}$
    \Repeat
        \State initialize $a$ ants on each UD;
        \For {each ant}
            \State construct a ring by repeatedly applying \eqref{prob}; 
            \State calculate the transmission time $T_{SR}^*$ 
            \Statex \quad \quad \, $\text{ of the ring;}$
            \If {$T_{SR}^* < T_{SC^*}^*$}
                \State Update the best ring, i.e. $r^*(\cdot) = r(\cdot)$;
            \EndIf
            \State $\triangle h_{i, j} += \frac{1}{T_{ring}^*}$;
        \EndFor
        \State update the pheromone by \eqref{update};
    \Until {reach the maximum iteration}
  \end{algorithmic}  
  \label{algorithm1}
\end{algorithm}  

In scenarios where the physical movement of UDs is not negligible, we may have different optimal rings in different communication rounds. At the beginning of each round, the BS collects the locations of UDs. Then, the BS invokes the ACO algorithm to determine the optimal logical ring and send it back to UDs. Such exchange of location information can be completed almost instantly as the data involved are very small, and thus would not have a significant impact on the implementation of our proposed policy.

Meanwhile, it is reasonable to consider that the UDs are static within a single communication round. The complexity of the ACO algorithm is polynomial, or in specific, $O(aN^2t)$, which makes the algorithm scalable for relatively large networks. Our numerical results also show that the average period of a single communication round is short. Therefore, the location of UDs and channel condition for any transmission pairs can be regarded as constant.

\section{Experimental Results}

\begin{figure*}[htbp]
    \centering
    \hspace{-0.1cm}
    \subfigure[\vspace{-0.05cm}\textit{The average total communication time in each round obtained by different communication schemes.}]{
    \includegraphics[scale=0.38]{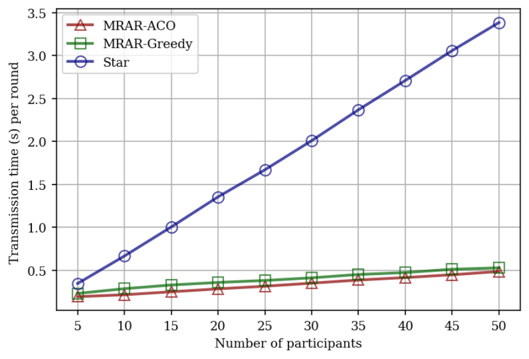}
    \label{line}}
    \subfigure[\textit{The average total communication time in each round with probability of D2D link failures.}]{\includegraphics[scale=0.38]{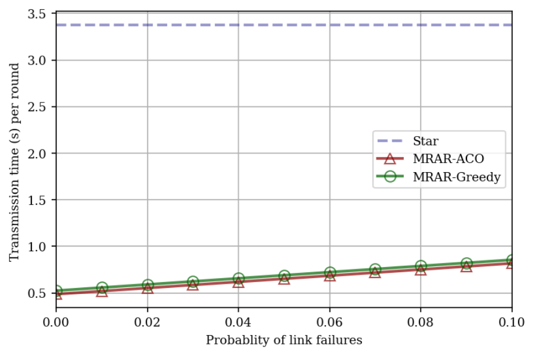}
    \label{inter}}
    \subfigure[\textit{Performance comparison of FedAvg with star and MRAR-based topologies, and DSGD with different preset distances.}]{\includegraphics[scale=0.222]{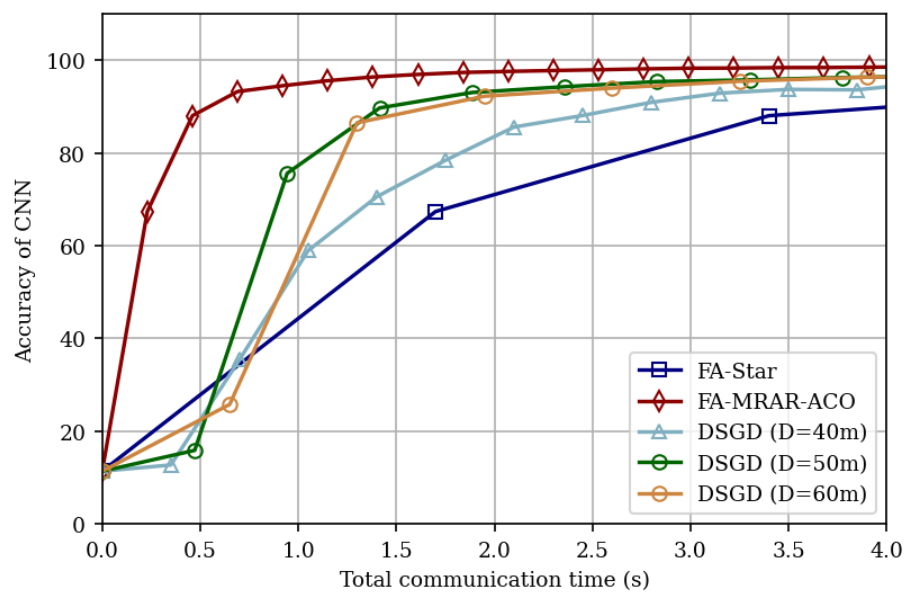}\label{DSGD}}
    \vspace{-0.05cm}
    \caption{\textit{Performance comparison of different communication schemes.}}
\end{figure*}
\subsection{Simulation setup}
We consider a $400m \times 400m$ squared area with a single BS at the central of the area.
The positions of UDs are uniformly and independently distributed in the area. 
We consider that all UDs have the same transmission power $p = 0.1$W. Other relevant network parameters are set as $M = 10$Mb, $B = 100$MHz, $\alpha = 4$ and $N_0 = -90$dBm. 

We consider two baseline policies referred to as \emph{star} and \emph{greedy}, respectively. In star, we consider the previously described traditional communication approach where all UDs only communicate with the BS to exchange the model parameters. In greedy, we consider a ring topology where the connection scheme is constructed by repeatedly connecting to the closest UD that has not been visited. We here reiterate that our approach focuses on the communication scheme of FL in wireless networks and does not alter the training algorithm (e.g. FedSGD, FedAvg, FedProx, FedPAQ) itself. 
Therefore, we can compare the performances of our proposed scheme and the baselines only by the transmission time of each round. Besides, we apply the optimal bandwidth allocation policy to every scenario to guarantee a more fair comparison. In the previous section, we have intuitively explained that the transmission time will not be longer than any other bandwidth allocation methods (e.g., all UDs are evenly allocated bandwidth) under all schemes. 

For the ACO, we set $\beta = 2$, $\gamma = 2$ and $\rho = 0.8$. We initialize $a = 10$ ants for each UD in each iteration and perform $t = 30$ iterations in total. 

\subsection{Numerical results}
We use $T_{\text{uStar}}^*$ as an optimistic estimation for 
$T_{\text{star}}^*$ in the following results. We simulate 50 cases for each scenario with the same number of UDs, and comparisons of averaging transmission time per round are demonstrated in Fig.~\ref{line}. We can conclude from the results that the MRAR architecture constructed by ACO (MRAR-ACO) outperforms the other two schemes in almost all scenarios. Meanwhile, the gaps between the transmission time of star and MRAR schemes (including greedy and ACO) become larger as the number of UDs increases. This is consistent with our previous analysis that, under MRAR, the total transmission time increases sublinearly with the density of UDs, which is superior to the linear increase in star. 

To demonstrate that our proposed MRAR is robust to link failures, we fix $K = 50$ and adjust the link failure probability. We simulate $50$ cases for each data point in Fig.~\ref{inter}. While the average transmission time per round by MRAR increases linearly with the failure probability, it is still significantly less than Star in all cases. This shows that our approach is resilient to D2D link failures and is still applicable to situations where D2D communication is relatively unstable.

In addition, we compare the performance between MRAR and DSGD, an existing approach we mentioned earlier that aggregates locally trained parameters in a different decentralized manner~\cite{xing2021federated}. The simulation results shown in Fig.~\ref{DSGD} are based on a trained CNN model by MINIST with 1,181,792 parameters. 
In DSGD, each UD broadcasts its local model to other UDs within a preset distance $D$, of which the value would affect the communication efficiency. For more intuitive comparisons,
we present the training accuracy versus total communication time of DSGD with different values of $D$, together with FedAvg implemented by our proposed MRAR formed by ACO (FA-MRAR-ACO), and FedAvg under the traditional star topology (FA-Star). As shown in Fig.~\ref{DSGD}, a larger value of $D$ in DSGD generally leads to better communication efficiency, as reflected by higher levels of final training accuracy and a faster convergence rate. The reason is that, while a shorter preset broadcast distance can shorten the communication time per round in DSGD, more rounds are required until convergence as a smaller $D$ naturally results in a larger number of smaller local clusters. On the other hand, while DSGD can improve communication efficiency compared to FA-Star under certain circumstances, our proposed FA-MRAR-ACO outperforms DSGD with all values of $D$. That is, FA-MRAR-ACO can achieve a very high level of training accuracy in a relatively short period, by both effectively reducing the communication time per round, and controlling the number of rounds required to achieve convergence.

\section{Conclusion}
The paper investigated the communication efficiency for FL in OFDMA wireless networks. We showed that the traditional star communication scheme for FL was inefficient especially when UDs were densely distributed. To overcome the problem, we proposed the MRAR architecture for FL parameter aggregation, and applied the ACO to optimize the formation of the ring. Analytical proof and simulation results showed that our proposed approach achieved significant improvement in both efficiency and robustness compared with baseline policies. 



%% file: main.bbl
\begin{thebibliography}{1}
\providecommand{\url}[1]{#1}
\csname url@samestyle\endcsname
\providecommand{\newblock}{\relax}
\providecommand{\bibinfo}[2]{#2}
\providecommand{\BIBentrySTDinterwordspacing}{\spaceskip=0pt\relax}
\providecommand{\BIBentryALTinterwordstretchfactor}{4}
\providecommand{\BIBentryALTinterwordspacing}{\spaceskip=\fontdimen2\font plus
\BIBentryALTinterwordstretchfactor\fontdimen3\font minus
  \fontdimen4\font\relax}
\providecommand{\BIBforeignlanguage}[2]{{%
\expandafter\ifx\csname l@#1\endcsname\relax
\typeout{** WARNING: IEEEtran.bst: No hyphenation pattern has been}%
\typeout{** loaded for the language `#1'. Using the pattern for}%
\typeout{** the default language instead.}%
\else
\language=\csname l@#1\endcsname
\fi
#2}}
\providecommand{\BIBdecl}{\relax}
\BIBdecl

\bibitem{li2020federated}
T.~Li, A.~K. Sahu, A.~Talwalkar, and V.~Smith, ``{Federated learning:
  Challenges, methods, and future directions},'' \emph{IEEE Signal Processing
  Mag.}, vol.~37, no.~3, pp. 50--60, 2020.

\bibitem{8917724}
J.~Mills, J.~Hu, and G.~Min, ``Communication-efficient federated learning for
  wireless edge intelligence in {IoT},'' \emph{IEEE Internet of Things
  Journal}, vol.~7, no.~7, pp. 5986--5994, 2020.

\bibitem{9177084}
A.~M. Elbir and S.~Coleri, ``{Federated Learning for Hybrid Beamforming in
  mm-Wave Massive MIMO},'' \emph{IEEE Commun. Lett.}, vol.~24, no.~12, pp.
  2795--2799, 2020.

\bibitem{liu2020client}
L.~Liu, J.~Zhang, S.~Song, and K.~B. Letaief, ``Client-edge-cloud hierarchical
  federated learning,'' in \emph{IEEE ICC 2020}, pp. 1--6.

\bibitem{7051236}
Y.~Wu, J.~Wang, L.~Qian, and R.~Schober, ``{Optimal Power Control for Energy
  Efficient D2D Communication and Its Distributed Implementation},'' \emph{IEEE
  Commun. Lett.}, vol.~19, no.~5, pp. 815--818, 2015.

\bibitem{xing2021federated}
H.~Xing, O.~Simeone, and S.~Bi, ``{Federated Learning Over Wireless
  Device-to-Device Networks: Algorithms and Convergence Analysis},'' \emph{IEEE
  J. Sel. Areas Commun.}, vol.~39, no.~12, pp. 3723--3741, 2021.

\bibitem{9796785}
M.~Yu, Y.~Tian, B.~Ji, C.~Wu, H.~Rajan, and J.~Liu, ``{GADGET: Online Resource
  Optimization for Scheduling Ring-All-Reduce Learning Jobs},'' in \emph{Proc.
  of IEEE INFOCOM 2022 - IEEE Conference on Computer Communications}, 2022, pp.
  1569--1578.

\bibitem{reisizadeh2020fedpaq}
A.~Reisizadeh, A.~Mokhtari, H.~Hassani, A.~Jadbabaie, and R.~Pedarsani,
  ``{FedPAQ: A Communication-Efficient Federated Learning Method with Periodic
  Averaging and Quantization},'' in \emph{International Conference on
  Artificial Intelligence and Statistics}, 2020, pp. 2021--2031.

\bibitem{dorigo2006ant}
M.~Dorigo, M.~Birattari, and T.~Stutzle, ``{Ant Colony Optimization},''
  \emph{IEEE computational intelligence magazine}, vol.~1, no.~4, pp. 28--39,
  2006.

\end{thebibliography}
